\documentclass[12pt]{article}
\usepackage{mathrsfs}
\usepackage{amsfonts}
\usepackage{amssymb}
\usepackage{amsthm}
\usepackage{amsmath}
\usepackage{latexsym}
\usepackage{color}
\begin{document}
\newtheorem{lemma}{Lemma}
\newtheorem{pron}{Proposition}
\newtheorem{thm}{Theorem}
\newtheorem{Corol}{Corollary}
\newtheorem{exam}{Example}
\newtheorem{defin}{Definition}
\newtheorem{rem}{Remark}
\newcommand{\la}{\frac{1}{\lambda}}
\newcommand{\sectemul}{\arabic{section}}
\renewcommand{\theequation}{\sectemul.\arabic{equation}}
\renewcommand{\thepron}{\sectemul.\arabic{pron}}
\renewcommand{\thelemma}{\sectemul.\arabic{lemma}}
\renewcommand{\thethm}{\sectemul.\arabic{thm}}
\renewcommand{\theCorol}{\sectemul.\arabic{Corol}}
\renewcommand{\theexam}{\sectemul.\arabic{exam}}
\renewcommand{\thedefin}{\sectemul.\arabic{defin}}
\renewcommand{\therem}{\sectemul.\arabic{rem}}
\def\REF#1{\par\hangindent\parindent\indent\llap{#1\enspace}\ignorespaces}
\def\E{\mathrm E}
\def\d{\mathrm d}
\def\P{\mathrm P}
\def\R{\mathbb{R}}

\title{\large\bf A Breiman's theorem for
conditional dependent random vector and its applications to risk theory}
\author{\small Zhaolei Cui$^1\ $\thanks{E-mail: cuishaol@126.com}\ \ \ and\ \ \ Yuebao Wang$^2$
\thanks{Corresponding author. E-mail: ybwang@suda.edu.cn}\\
{\small\it 1. School of mathematics and statistics, Changshu Institute of Technology, Suzhou, 215000, China}\\
{\small\it 2. Department of Mathematics, Soochow University, Suzhou 215006, China}}
\date{}

\maketitle {\noindent\small {\bf Abstract }}\\

{\small In this paper, we give a Breiman's theorem for conditional dependent random vector,
where one component has a regularly-varying-tailed distribution with the index $\alpha\ge0$,
while the other component is non-negative 
with a more relaxed moment condition.
This result substantially extends and improves some existing related results, such as
Theorem 2.1 of Yang and Wang (Extremes,\ 2013). 
We also provide some concrete examples, some interesting properties and a construction method of conditional dependent random vector.
Finally, we apply the above Breiman's theorem to risk theory,
and obtain two asymptotic estimates of the finite-time ruin probability
and the infinite-time ruin probability of a discrete-time risk model,
in which the corresponding net loss and random discount are conditionally dependent.\\



\section{\normalsize\bf Introduction}

\setcounter{equation}{0}\setcounter{defin}{0}
\setcounter{thm}{0}\setcounter{Corol}{0}\setcounter{lemma}{0}\setcounter{pron}{0}

It is well known that Breiman's theorem is an important tool
to deal with the tail distribution of product of two components of a random vector,
and it plays a key role in many application fields, such as risk theory, see Section 3 below.

In this paper, unless otherwise stated, let $(X,Y)$ be a random vector
with the marginal distributions $F$ on $(-\infty,\infty)$ and $G$ on $(0,\infty)$ respectively.
Then $Z=XY$ represents the product of these two components and its distribution denoted by $H$.
In risk theory, $X$ and $Y$ can be interpreted as the net loss 
and random discount, respectively.
Both of them are the basic research objects of risk theory.
The Breiman's theorem studies the asymptotics relations of the tail distributions of $Z$ and $X$.
Therefore, we should first introduce some concepts and symbols related to distribution.

Let $V$ be a distribution, then we say that $\overline{V}=1-V$ is the tail distribution of $V$.
If $V$ is supported on $(-\infty,\infty)$ \big(including $[0,\infty)$ and $(0,\infty)$\big), then
$$\overline{V}(x)>0\ \ \ \ \text{for all}\ \ x\in(-\infty,\infty).$$
Recall that a distribution $V$ on $(-\infty,\infty)$ is said to be regularly-varying-tailed with the index $\alpha\ge0$, denoted by $V\in\mathscr{R}_{-\alpha}$, if $\overline{V}(xy)\sim{\overline V}(x)y^{-\alpha}$, that is
$$\lim\overline{V}(xy)\big/{\overline V}(x)=y^{-\alpha}\ \ \ \ \text{for each}\ y>0.$$
Hereafter, all limits refer to $x\to\infty$ unless otherwise stated.
For more details on regularly-varying-tailed distribution, see Bingham et al. \cite{BGT1987} etc.
For example,
$$V\in\mathscr{R}_{-\alpha}\ \ \ \text{if and only if}\ \ \ {\overline V}(x)\sim x^{-\alpha}L(x),$$
where $L(\cdot)$ is a positive slowly varying function at infinity satisfying
$$L(xy)\sim L(x)\ \ \ \ \ \ \text{for each}\ \ y>0.$$

The distribution class $\cup_{\alpha\ge0}\mathscr{R}_{-\alpha}$ is a proper subclass
of the following distribution class introduced by Chistyakov \cite{C1964}.
Say that $V$ on $[0,\infty)$ \big(including $(0,\infty)$\big) belongs to the subexponential distribution class $\mathscr{S}$, if $$\overline{V^{*2}}(x)\sim2\overline{V}(x),$$
where $V^{*2}$ is the convolution of $V$ with itself.
Say that $V$ on $(-\infty,\infty)$ belongs to the class $\mathscr{S}$, if $V^+\in\mathscr{S}$,
where $V^+(x)=V(x)\textbf{1}_{[0,\infty)}(x)$ for all $x\in(-\infty,\infty)$.

Another distribution class introduced by Feller \cite{F1971} also properly includes class $\cup_{\alpha\ge0}\mathscr{R}_{-\alpha}$.
Say that distribution $V$ on $(-\infty,\infty)$ belongs to the dominantly-varying-tailed distribution class $\mathscr{D}$,
if for each $t\in(0,1)$, $\overline{V}(tx)=O\big(\overline{V}(x)\big)$, that is
$$\limsup\overline{V}(tx)\big/\overline{V}(x)<\infty.$$
If $V\in\mathscr{D}$, then it has a useful property as follows, which plays a key role in many cases,
see, for example, Theorems 1.C, 2.A, 2.B,  \ref{thm201},  \ref{thm301} and  \ref{thm302} below.
\\
\\
\textbf{Proposition 1.A.}
A distribution $V\in\mathscr{D}$ if and only if for any distribution $U$ on $(-\infty,\infty)$
satisfying $\overline {U}(x)=o(\overline {V}(x))$, that is
$$\lim\overline {U}(x)\big/\overline {V}(x)=0,$$
there exists a positive function $g(\cdot)$ on $[0,\infty)$ such that
\begin{equation}\label{101}
g(x)\downarrow0,\ \ xg(x)\uparrow\infty\ \ {\rm and}\ \ \overline {U}\big(xg(x)\big)=o\big(\overline{V}(x)\big).
\end{equation}

In the proposition, the proof of the necessity is given by Lemma 3.3 of Tang \cite{T2008b},
and the proof of sufficiency is attributed to Proposition 3.1 of Zhou et al. \cite{ZWW2012}.

In addition, the class $\mathscr{S}$ is properly included in the following distribution class.
Say that distribution $V$ on $(-\infty,\infty)$ belongs to the long-tailed distribution class $\mathscr{L}$,
denoted by $V\in\mathscr{L}$, if for each $t\in(0,\infty)$,
$$\overline{V}(x-t)\sim\overline{V}(x).$$

Goldie \cite{G1978} pointed out that class $\mathscr{L}$ and class $\mathscr{D}$ cannot contain each other.
Further, for a distribution $V$, we denote a positive function class by
$$\mathscr{H}_V=\{h(\cdot):\ h(x)\uparrow\infty,\ \ x/g(x)\downarrow0\ \ {\rm and}\ \ \overline{V}(x-t)\big/\overline{V}(x)\to1
\ \ \text{uniformly for all}\ \ |t|\le h(x)\}.$$
Then the following research way for distributions is often used in many occasions.
\\
\\
\textbf{Proposition 1.B.}
(i) A distribution $V\in\mathscr{L}$ if and only if set $\mathscr{H}_V$ is not empty.

(ii) A distribution $V\in\mathscr{S}$ if and only if $V\in\mathscr{L}$ and for any $h(\cdot)\in\mathscr{H}_V$,
$$\int_{h(x)}^{x-h(x)}\overline{V}(x-y)V(dy)=o\big(\overline{V}(x)\big).$$

For the proofs of (i) and (ii), please refer to Lemma 2.5 of Cline and Samorodnitsky \cite{CS1994}
and Proposition 1 of Asmussen et al. \cite{AFK2003}, respectively.
For a systematic overview of subexponential distribution class and related distribution classes,
please refer to Embrechts et al. \cite{EKM1997}, Resnick \cite{R2007}, Borovkov and Borovkov \cite{BB2008} and Foss et al. \cite{FKZ2013}, and so on.
\\

Now we return to the main research objects and objectives of this paper.
We continue to use the above notation in the following text.

Many papers have been devoted to studying the tail asymptotics of the product $Z=XY$,
especially the asymptotic relationship between $\overline{H}(x)$ and $\overline{F}(x)$.

In the case that $X$ and $Y$ are independent of each other, the following result was first stated by Breiman \cite{B1965} with $\alpha\in[0,1]$.
And the complete conclusion can be found in Corollary 3.6. (iii) of Cline and Samorodnitsky \cite{CS1994}.
\\
\\
\textbf{Theorem 1.A.}
For the random vector $(X,Y)$, assume that $X$ and $Y$ are independent of each other.
If $F\in\mathscr{R}_{-\alpha}$ for some $\alpha\in[0,\infty)$
and $\E Y^{\alpha+\epsilon}<\infty$ for some $\epsilon>0$, then $H\in\mathscr{R}_{-\alpha}$ and
\begin{equation}\label{102}
\overline {H}(x)\sim\E Y^\alpha\overline {F}(x).
\end{equation}

Since then, some papers have refined this result, see, for example,
Jessen and Mikosch \cite{JM2006}, Denisov and Zwart \cite{DZ2007} and Resnick \cite{R2007}.
Among them, Proposition 2.1 of Denisov and Zwart \cite{DZ2007} weakens the moment condition of r.v. $Y$
under some restrictions for distributions $F$ and $G$.
\\
\\
\textbf{Theorem 1.B.}
For the random vector $(X,Y)$, assume that $X$ and $Y$ are independent of each other
and $F\in\mathscr{R}_{-\alpha}$ for some $\alpha\in[0,\infty)$
with a positive slowly varying function at infinity $L(\cdot)$.
If $\E Y^\alpha<\infty,\ {\overline G}(x)=o({\overline F}(x))$ and
\begin{equation}\label{103}
\lim\sup_{1\leq y\leq x}L(x/y)\big/L(x)<\infty,
\end{equation}
then $H\in\mathscr{R}_{-\alpha}$ and (\ref{102}) holds.
\\

When $F$ is not necessarily a regularly-varying-tailed,
some corresponding results can be found in Hashorva et al. \cite{HPT2010}
for $F$ belonging to the max-domain of attraction of the Gumbel or Weibull distribution
and for $G$ belonging to the max-domain of attraction of the Weibull distribution,
Arendarczyk and D\c{e}bicki \cite{AD2011} for Weibull distributions $F$ and $G$,
and Hashorva and Li \cite{HL2014} and Cui et al. \cite{COWW2018} for semi-regular-varying-tailed $F$, among others.

In addition, based on Embrechts and Goldie \cite{EG1980}, 
Cline and Samorodnitsky \cite{CS1994} and 
Tang \cite{T2008}, 
Xu et al. \cite{XCWC2018} gave an equivalent condition for $H\in\mathscr{S}$ under the premise that $F\in\mathscr{S}$.
In particular, if $G$ is supported by $[0,a]$ for some $a>0$,
then a related result can also refer to Cline and Samorodnitsky \cite{CS1994}.
Other qualitative results can be found in Tang \cite{T2006}, Liu and Tang \cite{LT2010} and so on.

In the case that $X$ and $Y$ are dependent, which is a common phenomenon in practice,
Maulik et al. \cite{MRR2002} extended Theorem 1.A to the case that $Y$ is asymptotically independent
of $X$ (in a sense stronger than the usual concept of asymptotic independence),
Jiang and Tang \cite{JT2011} extended it to the case that $(X,Y)$
follows the two-dimensional Farlie-Gumbel-Morgenstern (FGM) distribution,
see case 1 after Definition \ref{def2010} below.
Furthermore, Yang and Wang \cite{YW2013} considered the two-dimensional Sarmanov distribution,
see Sarmanov \cite{S1974}, which is more general than the FGM distribution.

\begin{defin}\label{def2010}
Say that the random vector $(X,Y)$ follows a bivariate Sarmanov distribution, if
\begin{equation}\label{104}
\P(X\leq x,Y\leq y)=\int_{-\infty}^x\int_{-\infty}^y\big(1+\theta\phi_1(u)\phi_2(v)\big)F(du)G(dv),\ \ \ \ \ x,y\in(-\infty,\infty),
\end{equation}
where $\phi_1(\cdot)$ and $\phi_2(\cdot)$ are two measurable functions,
called kernels of this distribution, and the $\theta$ is a real constant satisfying
\begin{equation}\label{105}
\E\phi_1(X)=\E\phi_2(Y)=0,
\end{equation}
and
\begin{equation}\label{106}
1+\theta\phi_1(x)\phi_2(y)\geq0\ \ \ \ {\rm for\ all}\ x\in D_X\ \ \text{and}\ \ y\in D_Y,
\end{equation}
where
$$D_X=\big\{x\in(-\infty,\infty):\ \P\big(X\in(x-\delta,x+\delta)\big)>0\ \ \ {\rm for\ each}\ \delta>0\big\}$$
and
$$D_Y=\big\{y\in(0,\infty):\ \P\big(Y\in(y-\delta,y+\delta)\big)>0\ \ \ {\rm for\ each}\ \delta>0\big\}.$$
\end{defin}
Clearly, if $\theta=0$, or $\phi_1(x)=0$ for all $x\in D_X$, or $\phi_2(y)=0$ for all $y\in D_Y$, then $X$ and $Y$ are independent.
So, we say that a random vector $(X,Y)$ follows a proper bivariate Sarmanov distribution,
if the parameter $\theta\neq0$, and the kernels $\phi_1(\cdot)$ and $\phi_2(\cdot)$
are not identical to $0$ in $D_X$ and $D_Y$, respectively.
For some more details on multivariate Sarmanov distributions, one can refer
to Lee \cite{L1996}, Kotz et al. \cite{KBJ2000}, among others.
Three common choices for the kernels $\phi_1(\cdot)$ and $\phi_2(\cdot)$ are listed as follows.

Case 1. $\phi_1(x)=1-2F(x)$ and $\phi_2(y)=1-2G(y)$ for all $x\in D_X$
and $y\in D_Y$, leading to the well-known FGM distribution;

Case 2. $\phi_1(x)=(e^{-x}-c_1){\bf1}_{[0,\infty)}(x)$ with $c_1=\E
e^{-X}{\bf1}_{\{X\geq0\}}\big/\P(X\geq0)$ and $\phi_2(y)=e^{-y}-\E
e^{-Y}$ for all $x\in D_X$ and $y\in D_Y$;

Case 3. $\phi_1(x)=x^p-EX^p$ and $\phi_2(y)=y^p-EY^p$ for all $x\in
D_X$ and $y\in D_Y$.\vspace{1mm}

Proposition 1.1 of Yang and Wang \cite{YW2013} remarks that the kernels are bounded for any proper bivariate Sarmanov distribution.
\\
\\
\textbf{Proposition 1.C.}
Assume that random vector $(X,Y)$ follows a proper bivariate Sarmanov distribution.
Then there exist two positive constants
$b_1$ and $b_2$ such that
$$|\phi_1(x)|\le b_1\ \ \text{for all}\ x\in D_X\ \ \text{and}\ \ |\phi_2(y)|\le b_2\ \ \text{for all}\ \ y\in D_Y.$$

Motivated by Theorem 1.B and Jiang and Tang \cite{JT2011}, using Proposition 1.A,
Theorem 2.1 of Yang and Wang \cite{YW2013} gives the following result.
\\
\\
\textbf{Theorem 1.C.}
Let $(X,Y)$ be a random vector with a bivariate Sarmanov distribution. 
Assume that $F\in\mathscr{R}_{-\alpha}$ for some $\alpha\ge0$
with a positive slowly varying function at infinity $L(\cdot)$ and the limit $\lim\phi_1(x)=d_1$ exists.
If either (i) $\E Y^{\alpha+\epsilon}<\infty$ for some $\epsilon>0$;
or (ii) $\E Y^\alpha<\infty,\ {\overline G}(x)=o\big({\overline F}(x)\big)$
and there exist a positive function $g(\cdot)$ such that (\ref{101}) holds and
\begin{equation}\label{107}
\lim\sup_{1\leq y\leq xg(x)}L(x/y)\big/L(x)<\infty,
\end{equation}
then it holds that
\begin{equation}\label{108}
{\overline H}(x)\sim\E \big(s(Y)Y^\alpha\big){\overline F}(x),
\end{equation}
where $s(y)=1+\theta d_1\phi_2(y)$ for all $y>0$. 
\\

Further, Remark 2.1 of Yang and Wang \cite{YW2013} shows that
both Cases 1 and Case 2 satisfy all the conditions of Theorem 1.C.
And Remark 2.2 of Yang and Wang \cite{YW2013} notes that condition (\ref{107}) is properly weaker than condition (\ref{103}).

We might as well call Theorems 1.A, 1.B and 1.C the Breiman's theorem.

Naturally, we prefer to get a Breiman's theorem in a general dependent structure between $X$ and $Y$.
In addition, we intuitively believe that there should be some relationship between moment of $Y$ and function $s(\cdot)$.
And we hope that this relationship may help to weaken the moment condition of $Y$ on some occasions.

In Section 2, we introduce a conditional dependent structure of random vector $(X,Y)$,
and we give some concrete examples of conditional dependent random vectors,
which show that the conditional dependent structure is relatively large.
We also get some interesting properties and a construction method of conditional dependent random vectors.
Then, using a simple way, we prove a new Breiman's theorem for conditional dependent random vector $(X,Y)$.
Finally, in Section 3, according to this Breiman's theorem,
we obtain the asymptotic estimations of finite-time ruin probability and infinite-time ruin probability
in a discrete-time risk model with conditional dependent net loss and random discount.

\section{\normalsize\bf Breiman's theorem for conditional dependent random vector}

\setcounter{equation}{0}\setcounter{defin}{0}\setcounter{rem}{0}\setcounter{exam}{0}
\setcounter{thm}{0}\setcounter{Corol}{0}\setcounter{lemma}{0}\setcounter{pron}{0}

\subsection{\normalsize\bf Concepts, Examples and Properties}

We firstly introduce the concept of conditional dependent structure of random vector,
for research and application of which, please see, for example,
Badescu et al. \cite{BCL2009}, Li et al. \cite{LTW2010}, Asimit and Badescu \cite{AB2010},  Jiang et al. \cite{JWCX2015}.

\begin{defin}\label{def201}
Let $(X,Y)$ be a random vector.
For each $x\in(-\infty,\infty)$ and $y\in D_Y$, 
say that
$$\P(X>x|Y=y)=\lim_{t\downarrow0}\frac{\P\big(X>x,Y\in[y,y+t)\big)}{\P\big(Y\in[y,y+t)\big)}
$$
is the conditional distribution of $X$ under condition $Y=y$, if the above limit exits.

Say that $(X,Y)$ is conditionally dependent (CD),
if there exists a positive measurable function $s(\cdot)$ on $[0,\infty)$ such that
\begin{eqnarray*}
\P(X>x|Y=y)\sim\overline{F}(x)s(y)
\end{eqnarray*}
uniformly for all $y\in D_Y$, that is
\begin{eqnarray}\label{201}
\lim\sup_{y\in D_Y}\big|\P(X>x|Y=y)\big/\big(\overline{F}(x)s(y)\big)-1\big|=0.
\end{eqnarray}
\end{defin}

\begin{rem}\label{rem2010}
(i). For each pair $x\in(-\infty,\infty)$ and $y\notin D_Y$, if we define
$$\P(X>x|Y=y)=\overline{F}(x)\ \ \ \ \text{and}\ \ \ \ s(y)=1,$$
then (\ref{201}) still holds. Therefore, we can also expand the definition from $y\in D_Y$ to  $y>0$.

(ii). Later, we give some concrete examples of CD random vectors,
which show that the CD dependent structure has a wide range.
At the same time, we discuss some interesting properties of CD random vectors as follow.

1. Is $s(\cdot)$ bounded? see Proposition \ref{pron204}.

2. Can $s(\cdot)$ weaken the moment condition of $Y$? see Example \ref{exam201}.

3. Whether the CD property of a random vector is related to the parameters of its joint distribution?
see Proposition \ref{203}.

4. Besides some specific CD random vectors, is there a method to construct CD random vectors?
see Proposition \ref{pron205}.

5. Is $s(y)=1$ equivalent to $X$ and $Y$ independent of each other? see Proposition \ref{pron206}.
\end{rem}

\begin{pron}\label{pron201}
Let $(X,Y)$ be a random vector with a bivariate Sarmanov joint distribution.
If $\lim\phi_1(x)=d_1$ holds for some $d_1\in(-\infty,\infty)$, $\phi_2(\cdot)$ is continuous on $(0,\infty)$ and

\begin{eqnarray}\label{202}
s(y)=1+\theta d_1\phi_2(y)>0\ \ \ \ \ \text{for all}\ y\in D_Y,
\end{eqnarray}
then $(X,Y)$ is CD with the above function $s(\cdot)$ on $D_Y$.
\end{pron}}

\proof For each $y\in D_Y$, because $\phi_2(\cdot)$ is continuous on $y$,
then for any $\varepsilon\in(0,1)$, there exists a $\delta=\delta\big(\phi_2(\cdot),y,\varepsilon\big)>0$ such that
$$|\phi_2(v)-\phi_2(y)|<\varepsilon\ \ \ \text{for any}\ |v-y|<\delta.$$
Let $t\in(0,\delta)$, then
\begin{eqnarray*}
&&\P\big(X>x,Y\in[y,y+t)\big)=\int_{x}^{\infty}\int_{[y,y+t)}\big(1+\theta\phi_1(u)\phi_2(y)\big)F(du)G(dv)\nonumber\\
&&\ \ \ \ \ \ \ \ \ \ \ \ \ \ \ \ \ \ \ \ \ \ \ \ \ \ \ \ \ \ \ \ \ \ \ \ \ \ \ \ +\int_{x}^{\infty}\int_{[y,y+t)}\theta\phi_1(u)(\phi_2(v)-\phi_2(y))F(du)G(dv)\nonumber\\
&<&(\text{or}>)\ \P\big(Y\in[y,y+t)\big)\Big(\int_{x}^{\infty}\big(1+\phi_1(u)\phi_2(y)\big)F(du)\nonumber\\
&&\ \ \ \ \ \ \ \ \ \ \ \ \ \ \ \ \ \ \ \ \ \ \ \ \ \ \ \ \ \ \ \ \ \ \ \ +\varepsilon\ (\text{or}-\varepsilon)\int_{x}^{\infty}|\theta\phi_1(u)|F(du)\Big).
\end{eqnarray*}
Further, according to Proposition 1.C, we have
\begin{eqnarray}\label{203}
&&\P(X>x|Y=y)=\lim_{t\downarrow 0}\P\big(X>x,Y\in[y,y+t)\big)\big/\P\big(Y\in[y,y+t)\big)\nonumber\\
&=&\int_{x}^{\infty}\big(1+\theta\phi_1(u)\phi_2(y)\big)F(du)\nonumber\\
&=&\int_{x}^{\infty}\big(1+\theta d_1\phi_2(y)\big)F(du)+\int_{x}^{\infty}\theta\phi_2(y)\big(\phi_1(u)-d_1\big)F(du)
\end{eqnarray}
Therefore, by (\ref{202}), (\ref{203}) and $\phi_1(x)\to d_1$, we have
\begin{eqnarray*}
&&\lim\sup_{y\in D_Y}\Big|\frac{\P(X>x|Y=y)}{\overline{F}(x)s(y)}-1\Big|
=\lim\sup_{y\in D_Y}\Big|\frac{\int_{x}^{\infty}\big(1+\theta d_1\phi_2(y)\big)F(du)}{\overline{F}(x)s(y)}\nonumber\\
&&\ \ \ \ \ \ \ \ \ \ \ \ \ \ \ \ \ \ \ \ \ \ \ \ \ \ \ \ \ \ \ \ \ \ \ \  \ \ \ \ \ \ \ \ \ +\frac{\int_{x}^{\infty}\theta\phi_2(y)\big(\phi_1(u)-d_1\big)F(du)}{\overline{F}(x)s(y)}-1\Big|\nonumber\\
&=&0,
\end{eqnarray*}
that is $(X,Y)$ is CD.
\hfill$\Box$\\

Now we pay attention to another dependent structure of two-dimensional random vector,
see, for example, Nelsen \cite{N2006}.

\begin{defin}\label{def202}
Say that a random vector $(X,Y)$ follows a bivariate Frank distribution, if for some $\theta>0$,
\begin{eqnarray}\label{204}
\P(X\leq x,Y\leq y)=-\frac{1}{\theta}\ln\bigg(1+\frac{(e^{-\theta F(x)}-1)(e^{-\theta G(y)}-1)}{(e^{-\theta}-1)}\bigg),
\ \ \ x,\ y\in(-\infty,\infty).
\end{eqnarray}
\end{defin}

\begin{pron}\label{pron202}
Let $(X,Y)$ be the random vector with a bivariate Frank joint distribution, 
then $(X,Y)$ is CD with the a positive measurable function
$s(\cdot)$ on $(0,\infty)$ defined in (\ref{205}) below.
\end{pron}

\proof We first show that the following positive measurable function $s(\cdot)$ is positive and bounded on $(0,\infty)$. In fact, because $\theta>0$,
\begin{eqnarray}\label{205}
0<\frac{\theta e^{-\theta}}{1-e^{-\theta}}\le s(y)=\frac{\theta e^{-\theta\overline{G}(y)}}{1-e^{-\theta}}\le\frac{\theta}{1-e^{-\theta}}<\infty\ \ \ \ \ \ \text{for all}\ y>0.
\end{eqnarray}

Secondly, we give the expression of $\P(X>x|Y=y)$.
From (\ref{204}) and
$$\lim_{v\to0}(e^{-v}-1)/v=-1,$$
for each pair $x\in(-\infty,\infty)$ and $y\in D_Y$, we have
\begin{eqnarray*}
&&\P(X>x|Y=y)=1-\lim_{t\downarrow 0}\frac{\P\big(X\leq x,Y\in[y,y+t]\big)}{\P\big(Y\in[y,y+t]\big)}\nonumber\\
&=&1+\lim_{t\downarrow 0}\frac{1}{\theta \P\big(Y\in[y,y+t]\big)}
\ln\Big(\frac{(e^{-\theta}-1)+(e^{-\theta F(x)}-1)(e^{-\theta G(y+t)}-1)}
{(e^{-\theta}-1)+(e^{-\theta F(x)}-1)(e^{-\theta G(y)}-1)}\Big)\nonumber\\
&=&1+\lim_{t\downarrow 0}\frac{1}{\theta P(Y\in[y,y+t])}\ln\Big(1+\frac{(e^{-\theta F(x)}-1)(e^{-\theta G(y+t)}-e^{-\theta G(y)})}
{(e^{-\theta}-1)+(e^{-\theta F(x)}-1)(e^{-\theta G(y)}-1)}\Big)\nonumber\\
&=&1+\lim_{t\downarrow 0}\frac{\big(e^{-\theta F(x)}-1\big)(e^{-\theta G(y+t)}-e^{-\theta G(y)})}
{\theta P(Y\in[y,y+t])\big((e^{-\theta}-1)+(e^{-\theta F(x)}-1)(e^{-\theta G(y)}-1)\big)}\nonumber\\
&=&1+\frac{(e^{-\theta F(x)}-1)e^{-\theta G(y)}}{(e^{-\theta}-1)+(e^{-\theta F(x)}-1)(e^{-\theta G(y)}-1)}
\lim_{t\downarrow 0}\frac{e^{-\theta\big(G(y+t)-G(y)\big)}-1}{\theta\big(G(y+t)-G(y)\big)}\nonumber\\
&=&1-\frac{(e^{-\theta F(x)}-1)e^{-\theta G(y)}}{(e^{-\theta}-1)+(e^{-\theta F(x)}-1)\big(e^{-\theta G(y)}-1\big)}\nonumber\\
&=&\frac{e^{-\theta}-e^{-\theta F(x)}}{(e^{-\theta}-1)+(e^{-\theta F(x)}-1)(e^{-\theta G(y)}-1)}.
\end{eqnarray*}

Finally, we prove that $(X,Y)$ is CD. By
$$\lim \big(e^{-\theta}-e^{-\theta F(x)})\big/\overline{F}(x)=-\theta e^{-\theta}$$
and
$$\lim e^{-\theta F(x)}=e^{-\theta},$$
we know that, for any $\varepsilon>0$, there exists $x_0=x_0(F,\theta,\varepsilon)>0$ such that, when $x\geq x_0$.
$$|(e^{-\theta}-e^{-\theta F(x)})\big/\theta e^{-\theta}\overline{F}(x)+1|<\varepsilon$$
and
$$e^{-\theta F(x)}-e^{-\theta}<\varepsilon(1-e^{-\theta})\big/(e^\theta-1).$$
Thus, for each $x\geq x_0$ and the above $\varepsilon$, we have

\begin{eqnarray*}
&&\sup_{y\in D_Y}\Big|\frac{\P(X>x|Y=y)}{\overline{F}(x)s(y)}-1\Big|
=\sup_{y\in D_Y}\Big|\frac{\P(X>x|Y=y)-\overline{F}(x)s(y)}{\overline{F}(x)s(y)}\Big|\nonumber\\
&\leq&\sup_{y\in D_Y}\big|\P(X>x|Y=y)-\overline{F}(x)s(y)\big|(1-e^{-\theta})\big/\big(\theta e^{-\theta}\overline{F}(x)\big)\nonumber\\
&=&\frac{e^{\theta}-1}{\theta \overline{F}(x)}\sup_{y\in D_Y}\Big|\frac{e^{-\theta}-e^{-\theta F(x)}}
{(e^{-\theta}-1)+(e^{-\theta F(x)}-1)(e^{-\theta G(y)}-1)}
-\frac{\overline{F}(x)\theta e^{-\theta \overline{G}(y)}}{1-e^{-\theta}}\Big|\nonumber\\
&=&\frac{e^{\theta}-1}{\theta \overline{F}(x)}\sup_{y>0}\Big|\frac{e^{-\theta}-e^{-\theta F(x)}
+\theta e^{-\theta}\overline{F}(x)-\theta e^{-\theta}\overline{F}(x)}{(e^{-\theta}-1)
+(e^{-\theta F(x)}-1)(e^{-\theta G(y)}-1)}-\frac{\overline{F}(x)\theta e^{-\theta \overline{G}(y)}}{1-e^{-\theta}}\Big|\nonumber\\
&\leq&(1-e^{-\theta})\Big|\frac{e^{-\theta}-e^{-\theta F(x)}}
{\theta e^{-\theta}\overline{F}(x)}+1\Big|\sup_{y\in D_Y}\Big|\frac{1}{(e^{-\theta}-1)
+(e^{-\theta F(x)}-1)(e^{-\theta G(y)}-1)}\Big|\nonumber\\
&&\ \ \ \ \ \ \ \ \ \ \ \ \ \ +(1-e^{-\theta})\sup_{y\in D_Y}\Big|\frac{1}{(e^{-\theta}-1)
+(e^{-\theta F(x)}-1)(e^{-\theta G(y)}-1)}+\frac{e^{\theta G(y)}}{1-e^{-\theta}}\Big|\nonumber\\
&<&\varepsilon+\sup_{y\in D_Y}\big|(e^{-\theta F(x)}-e^{-\theta})(e^{\theta G(y)}-1)\big|\big/(1-e^{-\theta})\nonumber\\
&\le&\varepsilon+\sup_{y\in D_Y}\big|(e^{-\theta F(x)}-e^{-\theta})\big|(e^{\theta}-1)\big/(1-e^{-\theta})\nonumber\\
&<&2\varepsilon,
\end{eqnarray*}
which implies $(X,Y)$ is CD with the above $s(\cdot)$ on $(0,\infty)$
by the arbitrariness of $\varepsilon$.
\hfill$\Box$\\

The third two-dimensional joint distribution is as follows, see, for example, Li et.al\cite{LTW2010}.

\begin{defin}\label{pron203}
Say that a random vector $(X,Y)$ follows a bivariate Ali-Mikhail-Haq (AMH) distribution, if for some $\theta\in[-1,1)$,

\begin{eqnarray}\label{206}
\P(X\leq x,Y\leq y)=F(x)G(y)\big/\big(1-\theta\overline{F}(x)\overline{G}(x)\big),
\ \ \ \ \ x,\ y\in(-\infty,\infty).
\end{eqnarray}
\end{defin}

\begin{pron}\label{203}
Let $(X,Y)$ be the random vector with a bivariate AMH joint distribution. 
When $\theta\in(-1,1)$, $(X,Y)$ is CD with the a positive measurable function
\begin{eqnarray}\label{207}
s(y)=1+\theta\big(1-2\overline{G}(y)\big)\ \ \ \ \ \text{for all}\ \ y>0.
\end{eqnarray}
When $\theta=-1$,
\begin{eqnarray}\label{208}
\P(X>x|Y=y)\sim\overline{F}(x)s(y)\ \ \ \text{for each}\ \ y>0
\end{eqnarray}
and
\begin{eqnarray}\label{209}
\lim\sup_{y\in D_Y\cap(0,y_0]}\big|\P(X>x|Y=y)\big/\big(\overline{F}(x)s(y)\big)-1\big|=0,
\end{eqnarray}
where $s(y)=2\overline{G}(y)$ for all $y>0$ and
$$y_0=\sup\{y>0:3\overline{G}^2(y)\ge 1\}.$$
However, $(X,Y)$ is not CD.
\end{pron}

\proof Along the proof line of Proposition \ref{pron202},
for each $\theta\in[-1,1)$, $x\in(-\infty,\infty)$ and $y\in D_Y$, we first get that
\begin{eqnarray}\label{210}
&&\P(X>x|Y=y)=1-\lim_{t\downarrow 0}\P\big(X\leq x,Y\in[y,y+t]\big)\big/\P\big(Y\in[y,y+t]\big)\nonumber\\
&=&\lim_{t\downarrow 0}\frac{1-\theta\overline{F}(x)\big(\overline{G}(y+t)+\overline{G}(y)\big)
+\theta^2\overline{F}^2(x)\overline{G}(y+t)\overline{G}(y)-F(x)+\theta F (x)\overline{F}(x)}{\big(1-\theta \overline{F}(x)\overline{G}(y+t)\big)\big(1-\theta \overline{F}(x)\overline{G}(y)\big)}\nonumber\\
&=&\frac{\big(1+\theta-2\theta \overline{G}(y)-\theta \overline{F}(x)+\theta^2\overline{F}(x)\overline{G}^2(y)\big)\overline{F}(x)}
{\big(1-\theta \overline{F}(x)\overline{G}(y)\big)^2}.
\end{eqnarray}
Further, for any $\varepsilon>0$, we take $x_0=x_0(F,\varepsilon)>0$ such that
$\overline{F}(x)\leq \varepsilon$ for all $x\geq x_0$.

Next, we prove that $(X,Y)$ is CD for $\theta\in(-1,1)$. In fact,
\begin{eqnarray*}
&&\sup_{y\in D_Y}\Big|\frac{\P(X>x|Y=y)}{\overline{F}(x)s(y)}-1\Big|\\
&=&\sup_{y\in D_Y}\Big|\frac{\overline{G}(y)\big(2-\theta \overline{F}(x)\overline{G}(y)\big)\big(1+\theta(1-2\overline{G}(y))\big)-\big(1-\theta \overline{G}^2(y)\big)}{\big(1+\theta\big(1-2\overline{G}(y))\big)\big(1-\theta \overline{F}(x)\overline{G}(y)\big)^2}\Big||\theta|\overline{F}(x)\\
&\leq& \sup_{y\in D_Y}\Big(\Big|\frac{\overline{G}(y)\big(2-\theta \overline{F}(x)\overline{G}(y)\big)}{\big(1-\theta \overline{F}(x)\overline{G}(y)\big)^2}\Big|+\Big|\frac{1-\theta \overline{G}^2(y)}{\big(1+\theta(1-2\overline{G}(y))\big)\big(1-\theta \overline{F}(x)\overline{G}(y)\big)^2}\Big|\Big)\overline{F}(x)\\
&\leq& \Big( \frac{2+|\theta|\varepsilon}{(1-|\theta|\varepsilon)^2}
+\frac{1+|\theta|}{(1-|\theta|)(1-|\theta|\varepsilon)^2}\Big)\varepsilon.
\end{eqnarray*}
Then we complete the proof by the arbitrariness of $\varepsilon$.

For $\theta=-1$, by (\ref{210}), (\ref{208}) holds clearly.
(\ref{209}) comes from the following fact.
For the above $s(\cdot)=2\overline{G}(\cdot)$, $\varepsilon$ and $x_0$, when $x\ge x_0$, we have
\begin{eqnarray*}
&&\sup_{y\in D_Y\cap(0,y_0]}\Big|\frac{\P(X>x|Y=y)}{\overline{F}(x)s(y)}-1\Big|
=\sup_{y\in D_Y\cap(0,y_0]}\Big|\frac{3\overline{G}^2(y)+2\overline{F}(x)\overline{G}^3(y)-1}
{2\overline{G}(y)\big(1+\overline{F}(x)\overline{G}(y)\big)^2}\Big|\overline{F}(x)\\
&\leq&\sup_{y\in D_Y\cap(0,y_0]}\Big(\big(3\overline{G}(y)+2\overline{F}(x)\overline{G}^2(y)\big)
\Big/\big(1+\overline{F}(x)\overline{G}(y)\big)^2\Big)\overline{F}(x)/2\\
&\leq&(3+2\varepsilon)\varepsilon/2.
\end{eqnarray*}
Finally, we prove $(X,Y)$ is not a CD. For the above $\varepsilon$ and $x_0$, we further set $\varepsilon<1/5$.
For each $x\ge x_0$,
there exists $y=y(F,G,\varepsilon)>y_0$ large enough such that
$$\overline{G}(y)<\overline{F}(x)(1-5\varepsilon)\big/\big(2\big(1+\varepsilon)^2\big),$$
which leads to that
$$\overline{G}(y)<\varepsilon\ \ \ \text{and}\ \ \ 1>3\overline{G}^2(y)+2\overline{F}(x)\overline{G}^3(y).$$
Therefore, for the above $x$ and $y$, we have
\begin{eqnarray*}\label{2110}
&&\Big|\frac{\P(X>x|Y=y)}{\overline{F}(x)s(y)}-1\Big|
=\frac{\big(1-3\overline{G}^2(y)-2\overline{F}(x)\overline{G}^3(y)\big)\overline{F}(x)}
{2\overline{G}(y)\big(1+\overline{F}(x)\overline{G}(y)\big)^2}\\
&\ge&(1-5\varepsilon)\overline{F}(x)\big/\big(2\big(1+\varepsilon)^2\overline{G}(y)\big)\\
&>&1,
\end{eqnarray*}
that is $(X,Y)$ is not a CD.
\hfill$\Box$\\

In Propositions \ref{pron201}, \ref{pron202} and \ref{203},  functions $s(\cdot)$ are bounded.
In fact, this is not an individual phenomenon.

\begin{pron}\label{pron204}
Let $(X,Y)$ be a CD random vector with a function $s(\cdot)$ in Definition \ref{def201}.
Then $s(\cdot)$ is bounded in $y\in D_Y$, that is there is a $C>0$ such that
$$s(y)\le C\ \ \ \ \ \ \ \text{for all}\ y\in D_Y.$$
\end{pron}

\proof First, according to Definition \ref{def201}, $\P(X>x|Y=y)\le1$ for all $y\in D_Y$.
And because $F$ on $(-\infty,\infty)$, $\overline{F}(x)>0$ for each $x\in(-\infty,\infty)$.
Thus, $s(\cdot)$ is finite by (\ref{201}).

Secondly, we prove $s(\cdot)$ is bounded in $D_Y$.
Otherwise, there exist $y_n\in D_Y,\ n\ge1$ such that
$$y_n\uparrow\infty \ \ \text{as}\ n\to\infty\ \ \text{and}\ \ s(y_n)>n\ \ \text{for all}\ \ n\ge1.$$
And by (\ref{201}), there exists $x_0=x_0(X,Y)>0$ such that, when $x\ge x_0$,
\begin{eqnarray*}
\overline{F}(x)s(y)<2P(X>x|Y=y)\le2\ \ \ \ \ \ \ \text{for all}\ y\in D_Y.
\end{eqnarray*}
Take an integer $n_0=n_0(F,G,x_0)$ large enough such that $\overline{F}(x_0)n_0>2.$ Then
$$\overline{F}(x_0)s(y_{n_0})>\overline{F}(x_0)n_0>2,$$
which contradicts the previous inequality with $x=x_0$ and $y=y_{n_0}$.
\hfill$\Box$\\

In addition to the above three concrete examples,
we also hope to find a more general way to construct some CD random vectors.

\begin{pron}\label{pron205}
Let $(\xi,\eta)$ be a random vector with two marginal distributions $U$ and $V$ on $(0,\infty)$.
Assume that $(\xi,\eta)$ is CD and the corresponding function $s_0(\cdot)$ on $(0,\infty)$ is continuous
such that it holds uniformly for each $y\in D_\eta$ that
$$P(\xi>x|\eta=y)\sim\overline{U}(x)s_0(y).$$
Further, let $X=\xi-\eta$ and $Y=\eta$.
If $U\in\mathscr{L}$, then random vector $(X,Y)$ still is CD with a positive and continuous function
$s(\cdot)$ on $(0,\infty)$ satisfying
$$s(y)=s_0(y)\big/\E s_0(Y)\ \ \ \ \ \ \text{for each}\ \ y\in D_Y=D_\eta.$$
\end{pron}

\proof Because $(\xi,\eta)$ is CD, by $U\in\mathscr{L}$,
according to the dominant convergence theorem, we have
\begin{eqnarray*}
\overline{F}(x)=\int_{D_\eta}\P\big(\xi>x+v|\eta=v)V(dv)
\sim\int_{D_\eta}\overline{U}(x+v)s_0(v)V(dv)\sim\overline{U}(x)\E s_0(\eta).
\end{eqnarray*}
Further, considering the continuity of $s_0(\cdot)$, for each $y\in D_\eta$ and $\varepsilon\in(0,1)$,
there exists a small $t_0=t_0\big(s_0(\cdot),y,\varepsilon\big)>0$ and a large $x_0=x_0(U,y,\varepsilon,t_0)>0$ such that,
when $0<t\le t_0$ and $x\ge x_0$,
\begin{eqnarray*}
&&\P\big(X>x,Y\in[y,y+t)\big)\le\int_{[y,y+t)}\P(\xi>x|\eta=v)V(dv)\\
&\le&(1+\varepsilon)\overline{U}(x)s_0(y)\P\big(\eta\in[y,y+t)\big)
\end{eqnarray*}
and
\begin{eqnarray*}
&&\P\big(X>x,Y\in[y,y+t)\big)\ge\int_{[y,y+t)}\P(\xi>x-t|\eta=v)V(dv)\\
&\ge&(1-\varepsilon)\overline{U}(x)s_0(y)\P\big(\eta\in[y,y+t)\big).
\end{eqnarray*}
Thus, for each $y\in D_\eta$, by the arbitrariness of $\varepsilon$ and $\E s_0(\eta)<\infty$, we have
\begin{eqnarray*}
&&\P(X>x|Y=y)=\lim_{t\downarrow0}\P\big(\xi-\eta>x,\eta\in[y,y+t)\big)\big/\P\big(\eta\in[y,y+t)\big)\\
&=&\overline{U}(x)s_0(y),\\
\end{eqnarray*}
that is $(X,Y)$ is CD with $s(\cdot)=s_0(\cdot)\big/\E s_0(\eta)$ on $D_Y=D_\eta$.
\hfill$\Box$\\

We have reason to believe that there are many ways to construct CD random vectors.
Therefore, we can say that there are many CD random vectors.
More importantly, however, we find a somewhat surprising conclusion from the above proposition.

\begin{pron}\label{pron206}
Let $(X,Y)$ be a CD random vector with the corresponding function $s(\cdot)$ on $(0,\infty)$.
If $X$ is independent of $Y$, then $s(y)=1$ for all $y\in D_Y$.
On the contrary, $s(y)=1$ for all $y\in D_Y$ cannot implies the independence of $X$ and $Y$.
\end{pron}
\proof We only need to prove the second conclusion.
In fact, in Proposition \ref{pron205}, if $\xi$ is independent of $\eta$,
then $s_0(y)=1$ for all $y\in D_\eta$. Thus, according the proposition,
we know that $s(y)=s_0(y)\E s_0(\eta)=1$ for all $y\in D_Y$.
However, $X=\xi-\eta$ is not independent of $Y=\eta$ clearly.
\hfill$\Box$\\

\subsection{\normalsize\bf A new Breiman's theorem}

In this subsection, we prove a new Breiman's theorem for CD random vector $(X,Y)$.

\begin{thm}\label{thm201}
Let $(X,Y)$ be CD random vector. Assume that marginal distributions $F\in\mathscr{R}_{-\alpha}$ for some $\alpha\ge0$ with a positive slowly varying function at infinity $L(\cdot)$ and $G$ satisfies condition that ${\overline G}(x)=o\big({\overline F}(x)\big)$.
If either (i) for some $\varepsilon>0$,
\begin{eqnarray}\label{211}
\E Y^{\alpha+\varepsilon}s(Y)<\infty;
\end{eqnarray}
or (ii) $L(\cdot)$ satisfies (\ref{107}) for some positive function $g(\cdot)$ in (\ref{101}) and
\begin{eqnarray}\label{212}
\E Y^\alpha s(Y)<\infty,
\end{eqnarray}
then
\begin{eqnarray}\label{213}
\overline{H}(x)=\P(XY>x)\sim \E Y^\alpha s(Y)\overline{F}(x).
\end{eqnarray}
\end{thm}

\proof We firstly prove (ii).

Because $\mathscr{R}_{-\alpha}\subset\mathcal{D}$ and ${\overline G}(x)=o({\overline F}(x))$,
there is a positive function $g(\cdot)$ satisfying (\ref{101}) according to Proposition 1.A.
We take $g(\cdot)$ to sprite
\begin{eqnarray}\label{214}
\overline{H}(x)=\Big(\int_0^{xg(x)}+\int_{xg(x)}^\infty\Big)P(X>x/y\mid Y=y)G(dy)=\P_1(x)+\P_2(x).
\end{eqnarray}

For $\P_2(x)$, by (\ref{101}) and $\overline{F}(x)=O\big(\overline{H}(x)\big)$, we have

\begin{eqnarray}\label{215}
\P_2(x)\leq\overline{G}\big(xg(x)\big)=o\big(\overline{F}(x)\big)=o\big(\overline{H}(x)\big).
\end{eqnarray}

For $\P_1(x)$, by (\ref{107}) and (\ref{101}),
there exist four positive constants $x_0=x_0(F)$ and
$$C_1=C_1(x_0)\le C_2=C_2(x_0)\le C_3=C_3(x_0),$$
when $x\ge x_0$, it holds uniformly for all $y\in(0,xg(x)]$ that
\begin{eqnarray}\label{216}
&&\P(X>x/y|Y=y)\textbf{1}_{(0,xg(x)]}(y)\le P(X>x|Y=y)\textbf{1}_{(0,1]}(y)\nonumber\\
&&\ \ \ \ \ \ \ \ \ \ \ \ \ \ \ \ \ \ \ \ \ \ \ \ \ \ \ \ \ \ \ \ \ \ \ \ \ \ \ \ \ \ \ \ \
+\P(X>x/y|Y=y)\textbf{1}_{(1,xg(x)]}(y)\nonumber\\
&\le&C_1\big(\overline{F}(x)s(y)\textbf{1}_{(0,1]}(y)+\overline{F}(x/y)s(y)\textbf{1}_{(1,xg(x)]}(y)\big)\nonumber\\
&\le&C_2\big(s(y)\textbf{1}_{(0,1]}(y)+y^\alpha s(y)L(x/y)\big/L(x)\textbf{1}_{(1,xg(x)]}(y)\big)\overline{F}(x)\nonumber\\
&\le&C_3\big(s(y)\textbf{1}_{(0,1]}(y)+y^\alpha s(y)\textbf{1}_{(1,\infty]}(y)\big)\overline{F}(x).
\end{eqnarray}
Thus, according to dominant convergence theorem, by (\ref{216}) and (\ref{212}), we know that
\begin{eqnarray}\label{217}
&&\lim\P_1(x)\big/\overline{F}(x)=\int_0^{\infty}\lim\big(\P(X_1>x/y|Y=y)\big/\overline{F}(x)\big)\textbf{1}_{(0,xg(x)]}(y)G(dy)\nonumber\\
&=&\int_0^{\infty}\lim\big(\overline{F}(x/y)\big/\overline{F}(x)\big)\textbf{1}_{(0,xg(x)]}(y)s(y)G(dy)\nonumber\\
&=&\E Y^\alpha s(Y).
\end{eqnarray}

Combining (\ref{214}), (\ref{215}) and (\ref{217}), we get (\ref{213}) immediately.
\\

Next, we prove (i). To this end, we only need to deal with $L(x/y)\big/L(x)$ in the third line of ({\ref{211}).
According to Potter's theorem, see, for example, Theorem 1.5.6 of Bingham et al. \cite{BGT1987},
for the  above $x_0$ large enough and the $\varepsilon$ in (\ref{206}), when $x\ge x_0$,
there is a $C_3=C_3(F,x_0,\varepsilon)\ge C_2$ such that
\begin{eqnarray}\label{218}
L(x/y)\big/L(x)\le C_3y^{\alpha+\varepsilon}.
\end{eqnarray}
Thus, by (\ref{218}) and (\ref{211}), according to dominant convergence theorem,
we know that (\ref{212}) still holds.

Combining (\ref{214}), (\ref{215}) and (\ref{217}), (\ref{213}) still is holds.
\hfill$\Box$\\

According to Theorem \ref{thm201} and Proposition \ref{pron201}, we can get the following results directly.

\begin{Corol}\label{Corol201}
Let $(X,Y)$ be a random vector with a bivariate Sarmanov joint distribution
defined by (\ref{104}) satisfying (\ref{105}) and (\ref{106}).
Assume that all conditions of Proposition \ref{pron201} are satisfied,
then (\ref{213}) holds with $s(\cdot)$ in (\ref{202}).
\end{Corol}

\begin{rem}\label{rem201}
(i) According to Proposition 1.C or Proposition \ref{pron204},
we know that, there are two positive constants $C_1$ and $C_2$ such that
$$\E Y^{\alpha+\varepsilon}s(Y)\le C_1\E Y^{\alpha+\varepsilon}\ \ \ \text{and}\ \ \ \E Y^{\alpha}s(Y)\le C_2\E Y^{\alpha},$$
that is the order of $\E Y^{\alpha+\varepsilon}s(Y)$ (or $\E Y^{\alpha}s(Y)$)
is not stronger than $\E Y^{\alpha+\varepsilon}$ (or $\E Y^{\alpha}$).
In the following, however, we give a concrete random vector $(X,Y)$,
which satisfies all the conditions of Corollary \ref{Corol201}
and the order of $\E Y^{\alpha}s(Y)$ is significantly lower than that of $\E Y^{\alpha}$.
In other words, this result can substantially improve the moment condition of Theorem 1.C for many casses.

\begin{exam}\label{exam201}
Let $(X,Y)$ be a random vector with a bivariate Sarmanov joint distribution defined by (\ref{104}).

In (\ref{104}), we take a positive measurable function $\phi_1(\cdot)$ on $(-\infty,\infty)$
such that $\E \phi_1(X)=0$ and $\phi_1(x)\to d_1>0$.
Further, if we take
$$\theta=1\big/\big(d_1\E(1+1/Y)\big)\ \ \ \text{and}\ \ \ \phi_2(y)=1\big/(1+y)-1\big/(\theta d_1)\ \ \ \ \ \text{for all}\ y>0,$$
then $\theta>0$, $\phi_2(\cdot)$ is a continuous function on $(0,\infty)$,
$$s(y)=1+\theta d_1\phi_2(y)=\theta d_1\big/(1+y)>0\ \ \ \ \ \ \text{for all}\ \ y>0$$
and
$$\E\phi_2(Y)=\E1\big/(1+Y)-1\big/(\theta d_1)=0.$$
Thus, (\ref{105}) and (\ref{106}) are satisfied. At this time,
$$\E Y^{\alpha}s(Y)=\E Y^{\alpha}\big/(1+Y)\le G(1)+\E Y^{\alpha-1}\textbf{\emph{1}}_{\{Y>1\}}.$$
\end{exam}

(ii) Theorem 1.C does not require conditions that $\phi_2(\cdot)$ to be continuous and $s(y)>0$ for all $y>0$.
Thus, both Corollary \ref{Corol201} and Theorem 1.C have their own independent values.
Of course, in our opinion, it is more substantial to reduce the order of moments of $Y$.
\end{rem}

And according to Theorem \ref{thm201}, Proposition \ref{pron202} and Proposition \ref{pron205},
by (\ref{205}), we also have the three following results.

\begin{Corol}\label{Corol202}
Let $(X,Y)$ be a random vector with a bivariate Frank joint distribution defined by (\ref{204}).
Assume that all conditions of Proposition \ref{pron202} are satisfied,
then (\ref{213}) holds with $s(\cdot)$ in (\ref{205}).
\end{Corol}

\begin{Corol}\label{Corol203}
Let $(X,Y)$ be a random vector with a bivariate AMH joint distribution defined by (\ref{206}).
Assume that all conditions of Proposition \ref{pron203} are satisfied,
then (\ref{213}) holds with $s(\cdot)$ in the proposition.
\end{Corol}

\begin{Corol}\label{Corol204}
Let $(X,Y)$ be a random vector defined by Proposition \ref{pron205}.
Assume that all conditions of the proposition are satisfied,
then (\ref{213}) holds with $s(\cdot)$ in the proposition.
\end{Corol}

\section{\normalsize\bf Applications to risk theory}
\setcounter{equation}{0}\setcounter{rem}{0}
\setcounter{thm}{0}\setcounter{Corol}{0}\setcounter{lemma}{0}\setcounter{pron}{0}

We first introduce a discrete-time risk model which was proposed by Nyrhinen \cite{N1999, N2001}.

In the above model, within period $i$, the net insurance loss (i.e. the total claim amount
minus the total premium income) is denoted by a r.v. $X_i$ with a common distribution $F$ on $(-\infty,\infty)$.
And the random discount factor $Y_i$ from time $i$ to time $i-1$ with a common distribution $G$ on $(0,\infty)$.
As in Norberg \cite{N1999}, $X_n,\ n\geq1$ and $Y_n,\ n\geq1$ are called by insurance risks and financial risks, respectively.
For each integer $n\ge1$, the sum
\begin{equation}\label{301}
S_n=\sum_{i=1}^nX_i\prod_{j=1}^iY_j
\end{equation}
represents the stochastic discount value of aggregate net losses up to time $n$.
Let $x\geq0$ be the initial wealth of the insurer, then the finite-time ruin probability at time $n$ and
the infinite-time ruin probability can be respectively defined by
\begin{eqnarray}\label{302}
&&\psi_n(x)=\P\big(\max_{1\leq k\leq n}S_k>x\big)
\end{eqnarray}
and
\begin{eqnarray}\label{303}
&&\psi(x)=\lim_{n\to\infty}\psi(x,n)=\P\big(\sup_{n\geq1}S_n>x\big).
\end{eqnarray}

On the study of the asymptotics of $\psi_n(x)$ and $\psi(x)$,
when $\{X_i:i\ge1\}$ and $\{Y_i:i\ge1\}$ are independent of each other,
please refer to Tang and Tsitsiashvili \cite{TT2003, TT2004}, Goovaerts et al. \cite{GKLT2005},
Wang et al. \cite{WSZ2005}, Wang and Tang \cite{WT2006}, Zhang et al. \cite{ZSW2009}, Chen
and Yuen \cite{CY2009}, Shen et al. \cite{SLZ2009}, Gao and Wang \cite{GW2010}, Tang et
al. \cite{TVY2011}, Yi et al. \cite{YCS2011}, Zou et al. \cite{ZWW2012}, Cheng et al. \cite{CNPW2012}
and Xu et al. \cite{XCWC2018} among others.

Let vector $(X,Y)$ be an independent copy of $(X_i,Y_i)$ for all $i\ge1$.
When $(X,Y)$ follows a certain dependent structure, there are also some corresponding results.
Among them, Chen \cite{C2011} has obtained the following result that, for each fixed $n$,
\begin{equation}\label{304}
\psi(x,n)\sim\sum_{i=1}^n\P\Big(X_i\prod_{j=1}^iY_j>x\Big)
\end{equation}
under some conditions, where $(X,Y)$ follows a common bivariate FGM distribution.
And Theorem 4.1 and Theorem 4.2 of Yang and Wang \cite{YW2013} give the following two results
for $(X,Y)$ following a common bivariate Sarmanov distribution.

Hereafter, we denote by $H_i$ the distribution function of $X_i\prod_{j=1}^iY_j,\ i\geq1$. Clearly, $H_1=H$.
\\
\\
\textbf{Theorem 3.A.}
In the described discrete-time risk model, assume that $\{(X,Y),(X_i,Y_i):$ $i\ge1\}$ are i.i.d. random vectors.
Under the conditions Theorem 1.C, for each fixed integer $n\geq1$, it holds that
\begin{equation}\label{305}
\psi(x,n)\sim\Big(\E \big(s(Y)Y^\alpha\big)\big(1-\E^n Y^\alpha\big)\big/\big(1-\E Y^\alpha\big)\Big){\overline F}(x),
\end{equation}
where $s(y)=1+\theta d_1\phi_2(y),\ y\in D_Y$. Among them, we naturally agree that
\begin{equation}\label{306}
(1-\E^n Y^\alpha)/(1-\E Y^\alpha)=n,\ \ \ \ \text{if}\ \ \E Y^\alpha=1.
\end{equation}

\textbf{Theorem 3.B.}
In the described discrete-time risk model, assume that $\{(X,Y),(X_i,Y_i):$ $i\ge1\}$ are i.i.d. random vectors.
Under the conditions (i) of Theorem 1.C, for each fixed integer $n\geq1$, it holds that
\begin{equation}\label{307}
\psi(x)\sim\Big(\E \big(s(Y)Y^\alpha\big)\big/\big(1-\E Y^\alpha\big)\Big){\overline F}(x),
\end{equation}
where $s(y)=1+\theta d_1\phi_2(y),\ y\in D_Y$. Furthermore, (\ref{305}) holds uniformly for all $n\geq1$.
\\

Clearly, these two results describe the concrete influence of Sarmanov dependent structure of $(X,Y)$ on ruin probability.

In this section, using Theorem \ref{thm201},
we respectively give two asymptotic estimates of ruin probability of the discrete-time risk model,
where the net loss and random discount are CD.
Here, we can also find the influence of the CD structure on the ruin probability.

\begin{thm}\label{thm301}
In the described discrete-time risk model, assume that $\{(X,Y),(X_i,Y_i)$, $i\ge1\}$ are i.i.d.random vectors.
Under the conditions of Theorem \ref{thm201},
(\ref{305}) holds with $s(\cdot)$ in definition \ref{def201} for each fixed integer $n\geq1$.
\end{thm}

\proof We deal with case (i) and case (ii) together.

To prove (\ref{305}), we only need respectively to prove that
\begin{equation}\label{308}
\overline{H_i}(x)=P(X_iY_i\cdots Y_1>x)\sim\E^{i-1} Y^\alpha\E Y^\alpha s(Y){\overline F}(x)
\end{equation}
and
\begin{equation}\label{309}
\psi(x,n)\sim\sum_{i=1}^n\overline{H_i}(x)
\end{equation}
for all integers $1\le i\le n$ and each $n\ge1$ by induction.

We first prove (\ref{308}). According to Theorem \ref{thm201}, (\ref{308}) holds for $i=1$ obviously, which implies $H_1\in\mathscr{R}_{-\alpha}$ with a slowly varying function at infinity
$$L_1(x)\sim \E Y^\alpha s(Y)L(x),$$
where $L(\cdot)$ is a slowly varying function at infinity of distribution $F$.


Assume that (\ref{308}) holds for each integer $i\ge1$ and $H_i\in\mathscr{R}_{-\alpha}$
with a slowly varying function at infinity $L_i(\cdot)$ satisfying
\begin{equation}\label{310}
L_i(x)\sim\E^{i-1} Y^\alpha\E Y^\alpha s(Y)L(x).
\end{equation}
Clearly, replacing $L(\cdot)$ with $L_i(\cdot)$, condition (\ref{107}) is still satisfied and
$${\overline G}(x)=o\big({\overline F}(x)\big)=o\big(\overline {H_i}(x)\big)\ \ \ \ \text{for each}\ i\ge1.$$
Recall that $(X_i,Y_1),\ i\ge1$ are i.i.d. random vectors.
Thus, according to Theorem \ref{thm201}, by induction assumption, it holds that
\begin{eqnarray*}
\overline{H_{i+1}}(x)=\P \big((X_{i+1}Y_{i+1}\cdots Y_2)Y_1>x\big)
\sim\E Y^\alpha\overline{H_i}(x)\sim\E^{(i+1)-1} Y^\alpha\E Y^\alpha s(Y){\overline F}(x),
\end{eqnarray*}
that is (\ref{308}) and (\ref{310}) hold for $i+1$ and $H_{i+1}\in\mathscr{R}_{-\alpha}$.

Secondly, we prove (\ref{309}). To this end, we note that for each $n\geq1$,
$$\P(S_n>x)\leq\psi(x,n)\leq\P \Big(\sum_{i=1}^nX_i^+\prod_{j=1}^iY_j>x\Big)$$
and $(X^+,Y)$ also is CD.
This indicates that if we can establish
$$\P(S_n>x)\sim\sum_{i=1}^n\overline{H_i}(x),$$
then the asymptotic formula
\begin{equation}\label{311}
\P\left(\sum_{i=1}^nX_i^+\prod_{j=1}^iY_j>x\right)
\sim\sum_{i=1}^n\overline{H_i}(x)
\end{equation}
should hold as well. Furthermore, since $(X,Y),\ (X_i,Y_i),\ i\geq1$
are i.i.d. random vectors, it holds that
$$S_n=\sum_{i=1}^nX_i\prod_{j=1}^iY_j\stackrel{\d}{=}\sum_{i=1}^nX_i\prod_{j=i}^nY_j=T_n,\ \ \ n\geq1,$$
where $\stackrel{\d}{=}$ stands for equality in distribution.
Therefore, in order to prove (\ref{309}), we only need to prove that
\begin{equation}\label{312}
\P(T_n>x)=\overline{F_{T_n}}(x)\sim\sum_{i=1}^n\overline{H_i}(x)\sim\Big(\big(1-\E^n Y^\alpha)\E Y^\alpha s(Y)\big/(1-\E Y^\alpha)\Big){\overline F}(x)
\end{equation}
by induction on $n$.

Clearly, according to Theorem \ref{thm201}, by $F\in\mathscr{R}_{-\alpha}$, we know that
(\ref{312}) holds for $n=1$ and $F_{T_1}=H_1\in\mathscr{R}_{-\alpha}$.

For the sake of brevity, we only prove (\ref{312}) for $n=2$ in detail.
By the method that has been used many times before, we have
\begin{eqnarray}\label{313}
&&\overline{F_{T_2}}(x)\sim\int_0^{xg(x)}\P (X_1Y_1+X_2>x/y|Y_2=y)G(dy)\nonumber\\
&=&\int_0^{xg(x)}\Big(\int_{-\infty}^{x/y}+\int_{x/y}^\infty \P(X_2>x/y-z|Y_2=y)H_1(dy)\Big)G(dy)\nonumber\\
&=&\P_1(x)+\P_2(x).
\end{eqnarray}

Now we deal with $\P_1(x)$. To this end, by any $h(\cdot)\in\mathscr{H}_{F}$, we further split
\begin{eqnarray}\label{314}
&&\P_1(x)=\int_0^{xg(x)}\Big(\int_{-\infty}^{-h(x/y)}+\int_{-h(x/y)}^{h(x/y)}+
\int_{h(x/y)}^{x/y-h(x/y)}+\int_{x/y-h(x/y)}^{x/y}\Big)\nonumber\\
&&\ \ \ \ \ \ \ \ \ \ \ \ \ \ \ \ \ \ \ \ \ \ \ \ \ \ \ \ \ \ \ \cdot\P(X_2>x/y-z|Y_2=y)H_1(dz)G(dy)\nonumber\\
&=&\P_{11}(x)+\P_{12}(x)+\P_{13}(x)+\P_{14}(x).
\end{eqnarray}
By $F\in\mathscr{R}_{-\alpha}$, we know that
\begin{eqnarray}\label{315}
&&\P_{11}(x)\lesssim\int_0^{xg(x)}\overline{F}(x/y)s(y)H_1(-h(x/y)G(dy)\nonumber\\
&\lesssim&H_1\big(-h\big(1/g(x)\big)\overline{F}(x)\int_0^{xg(x)}\big(y^\alpha s(y)L(x/y)\big/L(x)\big)G(dy)\nonumber\\
&=&o(\overline{F}(x)).
\end{eqnarray}
According to Proposition 1.B (i), by $F\in\mathscr{R}_{-\alpha}$, we have
\begin{eqnarray}\label{316}
&&\P_{12}(x)\sim\int_0^{xg(x)}s(y)\int_{-h(x/y)}^{h(x/y)}\P(X_2>x/y-z)H_1(dy)G(dy)\nonumber\\
&\sim&\int_0^{xg(x)}\overline{F}(x/y)s(y)G(dy)\nonumber\\
&=&\E Y^\alpha s(Y)\overline{F}(x).
\end{eqnarray}
And according to Proposition 1.B (ii), by (\ref{308}) with $i=1$ and $F\in\mathscr{R}_{-\alpha}$, using integration by parts, we have
\begin{eqnarray}\label{317}
&&\P_{13}(x)\sim\int_0^{xg(x)}s(y)\int_{h(x/y)}^{x/y-h(x/y)}\P(X_2>x/y-z)H_1(dz)G(dy)\nonumber\\
&\sim&\int_0^{xg(x)}s(y)\int_{h(x/y)}^{x/y-h(x/y)}\overline{H_1}(x/y-u)F(du)G(dy)\nonumber\\
&\sim&\int_0^{xg(x)}s(y)o\big(\overline{F}(x/y)\big)G(dy)\nonumber\\
&=&o\big(\overline{F}(x)\big).
\end{eqnarray}
Finally, by (\ref{308}) with $i=1$ and $F\in\mathscr{R}_{-\alpha}$ again, we get that
\begin{eqnarray}\label{318}
&&\P_{14}(x)\leq\int_0^{xg(x)}\P\big(X_1Y_1\in(x/y-h(x/y),x/y]\big)G(dy)\nonumber\\
&=&\int_0^{xg(x)}o\big(\overline{H_1}(x/y)\big)G(dy)\nonumber\\
&=&o\big(\overline{F}(x)\big).
\end{eqnarray}

Next, we deal with $\P_2(x)$. On the one hand, it is obvious to get that
\begin{equation}\label{319}
\P_2(x)\leq\E Y^\alpha \overline{H_1}(x)\sim\E Y^\alpha\E Y^\alpha s(Y)\overline{F}(x).
\end{equation}
On the other hand, according to Fatou Lemma, by (\ref{308}), for any $a>1$, we have

\begin{eqnarray}\label{320}
&&\liminf\P_2(x)\big/\overline{H_1}(x)\geq\int_{-\infty}^\infty\liminf\int_{-\infty}^\infty \P(X_2>x/y-z|Y_2=y)\nonumber\\
&&\ \ \ \ \ \ \ \ \ \ \ \ \ \ \ \ \ \ \ \ \ \ \ \ \ \ \ \ \ \ \ \ \ \ \ \ \ \ \cdot\textbf{1}_{\{z>ax/y\}}H_1(dz)\textbf{1}_{\{y\leq xg(x)\}}\big/\overline{H_1}(x)G(dy)\nonumber\\
&\geq&\int_{-\infty}^\infty\liminf\int_{-\infty}^\infty\P\big(X_2>-(a-1)x/y|Y_2=y\big)\nonumber\\
&&\ \ \ \ \ \ \ \ \ \ \ \ \ \ \ \ \ \ \ \ \ \ \ \ \ \ \ \ \ \ \ \ \ \ \ \ \ \
\cdot\textbf{1}_{\{z>ax/y\}}H_1(dz)\textbf{1}_{\{y\leq xg(x)\}}\big/\overline{H_1}(x)G(dy)\nonumber\\
&=&\int_{-\infty}^\infty\liminf \P\big(X_2>-(a-1)x/y|y_2=y\big)\overline{H_1}(ax/y)\big/\overline{H_1}(x)G(dy)\nonumber\\
&=&\int_{-\infty}^\infty\liminf\overline{H_1}(ax/y)\big/\overline{H_1}(x)G(dy)\nonumber\\
&=&a^{-\alpha}\E Y^\alpha s(Y).
\end{eqnarray}
Then let $a \downarrow 1$, combining with (\ref{319}), we know that
\begin{equation}\label{321}
\P_2(x)\sim\E Y^\alpha \E Y^\alpha s(Y)\overline{F}(x).
\end{equation}

Hence by (\ref{313})-(\ref{318}) and (\ref{321}), (\ref{312}) is holds for $n=2$.

This completes the proof of Theorem \ref{thm301}.\hfill$\Box$\vspace{2mm}

\begin{thm}\label{thm302}
In the described discrete-time risk model, assume that $\{(X,Y),(X_i,Y_i)$, $i\ge1\}$ are i.i.d.random vectors.
Under the conditions of Theorem \ref{thm201} (i),
(\ref{307}) holds with $s(\cdot)$ in Definition \ref{def201} for each fixed $n\geq1$.
\end{thm}

The proof of this theorem is completely similar to the proof of Theorem 4.2 of Yang and Wang \cite{YW2013},
and we have omitted its details.
On the contrary, the proof of Theorem \ref{thm301}
is substantially different from the proof of Theorem 1.C.


\begin{thebibliography}{99}

\bibitem{AD2011} Arendarczyk M. and D\c{e}bicki K., 2011.
Asymptotics of supremum distribution of a Gaussian process over a Weibullian time. Bernoulli, 17(1), 194--210.

\bibitem{AB2010} Asimit, A., Badescu, A.L., 2010. Extremes on the discounted aggregate claims in a
time dependent risk model. Scand. Actuar. J. 2, 93–104.

\bibitem{AFK2003} Asmussen S., Foss, S. and Korshunov, D., 2003. Asymptotics for sums of random variables
with local subexponential behavior. J. Theor. Probab., 16, 489-518.

\bibitem{BCL2009} Badescu, A.L., Cheung, E.C.K., Landriault, D., 2009. Dependent risk models with
bivariate phase-type distributions. J. Appl. Probab. 46, 113–131.

\bibitem{BGT1987} Bingham, N. H., Goldie, C. M. and Teugels, J. L., 1987.
Regular Variation. Cambridge University Press, Cambridge.

\bibitem{BB2008} Borovkov A.A. and Borovkov, K.A. Asymptotic Analysis of random walks. Cambridge University Press, Cambridge, 2008.

\bibitem{B1965} Breiman, L., 1965. On some limit theorems similar to the arc-sin law. Theory Probab. Appl.
10, 323--331.

\bibitem{C2011} Chen, Y., 2011. The finite-time ruin probability with dependent insurance and
financial risks. J. Appl. Probab. 48, 1035--1048.

\bibitem{CY2009} Chen, Y. and Yuen, K. C., 2009. Sums of pairwise
quasi-asymptotically independent random variables with consistent
variation. Stochastic Models 25, 76--89.


\bibitem{C1964} Chistyakov V.P., 1964. A theorem on sums of independent, positive random variables
and its application to branching processes. Probability theory and its applications, 9, 640-648.

\bibitem{CS1994} Cline, D. B. H. and Samorodnitsky, G., 1994. Subexponentiality of the product of independent
random variables. Stochastic Process. Appl. 49, 75--98.

\bibitem{COWW2018} Cui Z., Omey E., Wang W. and Wang Y., 2018. Asymptotics of convolution
with the semi-regular-variation tail and its application to risk. Extremes, 21, 509-532.

\bibitem{DZ2007} Denisov, D. and Zwart, B., 2007. On a theorem of Breiman and a class of
random difference equations. J. Appl. Probab. 44, 1031--1046.

\bibitem{EG1980} Embrechts P. and M. C. Goldie, 1980. On closure and factorization properties of subexponential tails.
J. Austral. Math. Soc. (Ser. A), 29, 243--256.

\bibitem{EKM1997} Embrechts P., Kl\"{u}ppelberg C. and Mikosch T. Modelling Extremal Events for Insurance
and Finance. Springer, 1997.

\bibitem{F1971} Feller W., 1971. An Introduction to Probability Theory and its Applications, 2nd edn. vol. 2, Wiley, New York.

\bibitem{G1978} Goldie. C. M., 1978. Subexponential distributions and dominated variation tails.
J. Appl. Probab., 15, 440--442.

\bibitem{falk} Falk M., H\"usler, J. and Reiss R.-D., 2010. Laws of Small Numbers:
Extremes and Rare Events. DMV Seminar 23,  Third edition,
Birkh\"auser, Basel.

\bibitem{FKZ2013} Foss S., Korshunov D. and Zachary S. An Introduction to Heavy-tailed and Subexponential Distributions.
Springer, Second Edition, 2013.

\bibitem{GW2010} Gao, Q. and Wang, Y., 2010. Randomly weighted sums with dominated
varying-tailed increments and application to risk theory. J. Korean
Stat. Society 39, 305--314.

\bibitem{GKLT2005} Goovaerts, M. J., Kaas, R., Laeven, R. J. A., Tang, Q. and Vernic,
R., 2005. The tail probability of discounted sums of Pareto-like
losses in insurance. Scand. Actuar. J. 6, 446--461.

\bibitem{HL2014} Hashorva, E. and Li, J., 2014. Asymptotics for a discrete-time risk model with the emphasis
on financial risk. Probab. Eng. Inform. Sci. 28(4), 573-588.

\bibitem{HPT2010} Hashorva, E., Pakes, A. G. and Tang, Q., 2010. Asymptotics of random
contractions. Insurance Math. Econom. 47, 405--414.

\bibitem{JM2006} Jessen, A. H. and Mikosch, T., 2006. Regularly varying functions. Institut Math\'{e}matique.
Publications. Nouvelle S\'{e}rie 80, 171--192.

\bibitem{JT2011} Jiang, J. and Tang, Q., 2011. The product of two dependent random variables with regularly varying
or rapidly varying tails. Statist. Probab. Lett. 81, 957--961.

\bibitem{JWCX2015} Jiang, T., Wang, Y., Chen, Y., Xu, H., 2015. Uniform asymptotic estimate for finite-time ruin probabilities of a
time-dependent bidimensional renewal model. Insurance: Mathematics and Economics 64, 45-53.

\bibitem{KBJ2000} Kotz, S., Balakrishnan, N. and Johnson, N. L., 2000. Continuous
Multivariate Distributions. Vol. 1: Models and Applications. Wiley.

\bibitem{L1996} Lee, M. T., 1996. Properties and applications of the Sarmanov family of
bivariate distributions. Comm. Statist. Theory Methods 25,
1207--1222.

\bibitem{LTW2010} Li, J., Tang, Q., Wu, R., 2010. Subexponential tails of discounted aggregate claims in
a time-dependent renewal risk model. Adv. Appl. Probab. 42, 1126–1146.

\bibitem{LT2010} Liu, Y. and Tang, Q., 2010. The subexponential product convolution of
two Weibulltype distributions.  J. Aust. Math. Soc. 89, 277--288.

\bibitem{MRR2002} Maulik, K., Resnick, S. and Rootz\'{e}n, H., 2002. Asymptotic independence and a network traffic
model. J. Appl. Probab. 39, 671--699.

\bibitem{N2006} Nelsen R. B, 2006. An Introduction to Copulas, Second edition. Springer, New York.

\bibitem{N1999} Norberg, R., 1999. Ruin problems with assets and liabilities
of difusion type. Stoch. Process. Appl. 81, 255-269.

\bibitem{N1999} Nyrhinen, H., 1999. On the ruin probabilities in a general economic
environment. Stoch. Process. Appl. 83, 319--330.

\bibitem{N2001} Nyrhinen, H., 2001. Finite and in.nite time ruin probabilities in a
stochastic economic environment. Stoch. Process. Appl. 92,
265--285.

\bibitem{R2007} Resnick, S. I., 2007. Heavy-Tail Phenomena: Probabilistic and
Statistical Modeling. Springer, New York.

\bibitem{S1974} Sarmanov, I. O., 1974. New forms of correlation relationships between positive
quantities applied in hydrology. In: Mathematical Models in Hydrology.
International Association of Hydrological Sciences, Paris, 104-109.

\bibitem{SLZ2009} Shen, X., Lin, Z. and Zhang, Y., 2009. Uniform estimate for
maximum of randomly weighted sums with applications to ruin theory.
Methodol. Comput. Appl. Probab. 11, 669--685.

\bibitem{T2006} Tang, Q., 2006. The subexponentiality of products revisited. Extremes 9, 231--241.

\bibitem{T2008} Tang, Q., 2008. From light tails to heavy tails through multiplier. Extremes 11, 379--391.

\bibitem{T2008b} Tang, Q., 2008. Insensitivity to negative dependence of asymptotic tail probabilities
of sums and maxima of sums. Stoch. Anal. Appl. 26, 435--450.

\bibitem{TT2003} Tang, Q. and Tsitsiashvili, G., 2003. Precise estimates for
the ruin probability in finite horizon in a discrete-time model with
heavy-tailed insurance and financial risks. Stoch. Process.\
Appl.\ 108, 299--325.

\bibitem{TT2004} Tang, Q. and Tsitsiashvili, G., 2004. Finite- and infinite-time ruin
probabilities in the presence of stochastic returns on investments.
Adv. Appl. Probab. 36, 1278--1299.

\bibitem{TVY2011} Tang, Q., Vernic, R. and Yuan, Z., 2011. The finite-time ruin probability in the presence of
dependent extremal insurance and financial risks. .

\bibitem{WSZ2005} Wang, D., Su, C. and Zeng, Y., 2005. Uniform estimate for maximum
of randomly weighted sums with applications to insurance risk
theory. Science in China: Series A 48, 1379--1394.

\bibitem{WT2006} Wang, D. and Tang, Q., 2006. Tail probabilities of randomly
weighted sums of random variables with dominated variation.
Stoch. Models 22, 253--272.

\bibitem{XCWC2018} Xu, H., Cheng, F., Wang Y. and Cheng D., 2018. A necessary and sufficient condition for the subexponentiality
of product convolution. Adv. Appl. Prob., 50, 57-73.

\bibitem{YW2013} Yang Y. and Wang Y., 2013. Tail behavior of the product of two dependent random variables
with applications to risk theory. Extremes, 16, 55–74.

\bibitem{YCS2011} Yi, L., Chen, Y. and Su, C., 2011. Approximation of the tail
probability of randomly weighted sums of dependent random variables
with dominated variation. J. Math. Anal. Appl. 376, 365--372.

\bibitem{ZSW2009} Zhang, Y., Shen, X. and Weng, C., 2009. Approximation of the tail
probability of randomly weighted sums and applications. Stoch.
Process. Appl. 119, 655--675.

\bibitem{ZWW2012} Zhou, M., Wang, K. and Wang, Y., 2012. Estimates for the finite-time ruin probability with insurance
and financial risks. Acta Math. Appl. Sin., Engl. Ser., 28(4), 795–806.

\end{thebibliography}
\end{document}